\documentclass[11pt,notitlepage,twoside]{article}
\usepackage{latexsym}
\usepackage{amsmath}
\usepackage{amsfonts}
\usepackage{amssymb}

\newcommand{\e }{\varepsilon }

\newcommand{\be}{\begin{equation}}
\newcommand{\ee}{\end{equation}}
\newcommand{\ba}{\begin{align*}}
\newcommand{\ea}{\end{align*}}

\newtheorem{theorem}{Theorem}[section]
\newtheorem{definition}{Definition}[section]
\newtheorem{lemma}{Lemma}[section]
\newtheorem{proposition}{Proposition}[section]

\newtheorem{corollary}{Corollary}[section]

\author{Zakaria Boucheche\\
{\footnotesize  Laboratory of Applied Mathematics and Harmonic
Analysis (LR11-ES52)}\\{\footnotesize  (Faculty of Sciences of
Gabes. Gabes University)}}

\title { \Large \textbf{Another proof of a Lions type existence result}}

\begin{document}

\date{ }

\maketitle{\footnotesize} \noindent {\bf Abstract.} This paper
concerns a nonlinear elliptic equation involving a critical Sobolev
growth and a lower-order term. Under a Lions's condition, we prove
the existence of at least one positive solution. Our approach
consists in constructing a relatively compact Palais--Smale sequence
for the associated variational problem.

\noindent{ \footnotesize {{\it \\\\2010 Mathematics Subject
Classification.} \,\, 35J60, 35A15, 35B33, 37C10, 54D45, 35A01.}
}
\section{Introduction and main results}
\def\theequation{1.\arabic{equation}}\makeatother
\setcounter{equation}{0} We study the following nonlinear elliptic
partial differential equation with zero Dirichlet boundary condition
\begin{equation}\label{problem}
\begin{aligned}
-\Delta& u=K(x)u^{q}+\mu u \quad \mbox{in}\,\,\Omega,\\{}& u>0 \quad
\mbox{in} \,\,\Omega,\quad u=0 \quad \mbox{on} \,\,\partial\Omega,
\end{aligned}
\end{equation}where $\Omega\subset \mathbb{R}^n,\,n\geq 3,$ is a bounded domain with a smooth boundary
$\partial\Omega,$ $K(x)$ is a  continuous function in
$\bar{\Omega},$ $q+1=\frac{2n}{n-2}$ is the critical exponent for
the embedding $H_0^1\bigl(\Omega\bigr)$ into
$L^{q+1}\bigl(\Omega\bigr)$ and $0< \mu\leq \mu_1(\Omega),$ where
$\mu_1(\Omega)$ denotes the first eigenvalue of $(-\Delta)$ in
$H_0^1\bigl(\Omega\bigr),$ \\

One motivation to study this equation comes from its resemblance to
the well known scalar curvature problem on an $n$-dimensional closed
Riemannian manifold $(M^n,\,g_0),\,n\geq 3,$ which consists to find
a new metric $g$ conformally equivalent to $g_0$ with prescribed
scalar curvature $K(x)$ on $M^n;$ see, e.g., \cite{Au}.\\

Before setting forth the main existence result, let us introduce
some notations. Let $\langle\,,\,\rangle$ denotes the scalar product
defined on $H_0^1\bigl(\Omega\bigr)$ by $$\langle
u,v\rangle=\int_\Omega\nabla u\cdot \nabla v$$ and let $\|\cdot\|$
denotes its associated norm. Let $K_\infty,$ $S$ and $L_{K,\mu}$
denote the following constants \begin{eqnarray}
&K_\infty:=\sup_{\bar{\Omega}}(K),\quad S:=\text{inf}\{\|u\|^2,\quad
u\in H_0^1\bigl(\Omega\bigr)\,\,\text{and\,\,}\|u\|_{q+1}=1\},
\nonumber\\&L_{K,\mu}:=\text{inf}\{\|u\|^2-\mu \|u\|_2^2,\quad u\in
H_0^1\bigl(\Omega\bigr)\,\,\text{and\,\,}J(u)=1\},\label{inf}
\end{eqnarray}
where $J(u):=\int_\Omega K(x)|u(x)|^{q+1}\,\text{d}x\,\,
\text{and\,\,}\|u\|_p^p=\int_\Omega |u(x)|^p\,\text{d}x,$ for any
$p> 1.$ $S$ is known as the best Sobolev constant.\\ For our present
problem, we read the Lions type theorem as follows:
\begin{theorem}\label{Lions}{\it\,\,\,Let $n\geq 3.$ Assume that $K_\infty>
0$ and that $0< \mu< \mu_1(\Omega).$ If
\begin{equation}\label{zakaria}
L_{K,\mu}<\frac{1}{\bigl(K_\infty\bigr)^{\frac{n-2}{n}}}S,
\end{equation}then the problem \eqref{problem} has a
solution $u$ satisfying
$$[J(u)]^{\frac{n-2}{n}}L_{K,\mu}\leq\|u\|^2-\mu
\|u\|_2^2<\frac{1}{\bigl(K_\infty\bigr)^{\frac{n-2}{n}}}S[J(u)]^{\frac{n-2}{n}}.$$}\end{theorem}In
\cite{Li}, Lions introduced a concentration-compactness method,
which enabled him, from others, to study the loss of compactness
related to the constrained minimization
problem\begin{equation}\label{minimizing}\text{inf}\{\|u\|^2-\mu
\|u\|_2^2,\quad u\in
H_0^1\bigl(\Omega\bigr)\,\,\text{and\,\,}J(u)=1\}.\end{equation} In
\cite[Corollary 4.1]{Li1}, the author proved that the hypotheses of
Theorem \ref{Lions} are sufficient to ensure that any minimizing
sequence of this problem is relatively compact, and then a solution
to problem \eqref{problem} is regained at the level set
$\varpi:=(1/n)L_{K,\mu}^{n/2}$ for the functional $I_{K,\mu}$
defined in \eqref{PS0} below. See the proof of
\cite[Theorem I.2]{Li} for more details.\\

The problem \eqref{problem} enjoys a variational structure. Indeed
solution of \eqref{problem} corresponds to positive critical point
of the functional $I_{K,\mu}$ defined on $H_0^1(\Omega)$
by\begin{equation}\label{PS0} I_{K,\mu}(u)=\frac{1}{2}\int_{\Omega}
|\nabla u|^2 - \frac{1}{q+1} \int_{\Omega} K|u|^{q+1}- \frac{\mu}{2}
\int_{\Omega} u^2,\quad\forall\,\,u\in H_0^1(\Omega).
\end{equation}To resolve \eqref{problem} one
can think to the Palais--Smale (P-S for short) condition for
$I_{K,\mu}.$ Let $\partial I_{K,\mu}$ denotes the gradient of
$I_{K,\mu}.$
\begin{definition}\label{defPS}\,\,Let $c\in \mathbb{R}.$\\
1) Let $(u_k)_k$ be a sequence in $H_0^1(\Omega).$ We say that
$(u_k)_k$ is a P-S sequence at $c$ for $I_{K,\mu}$ if, up a
subsequence, $I_{K,\mu}\bigl(u_k\bigr)\rightarrow c$ and $\partial
I_{K,\mu}\bigl(u_k\bigr)\rightarrow 0$ strongly in $H^{-1}(\Omega).$\\
2) We say that $I_{K,\mu}$ satisfies the P-S condition at $c$ if any
P-S sequence at $c$ is relatively compact.
\end{definition}
As well-known in variational problems with critical exponent,
concentration phenomena can occur and violate the P-S condition at
some levels. Thus a local analysis becomes useful. In this
direction, a mountain-pass procedure was introduced in \cite{BN}: In
order to obtain a P-S sequence, the authors used, as a key tool, an
Ambrosetti--Rabinowitz type result \cite[Theorem 2.2]{BN} . For an
adaptation of this procedure to the present problem, we refer the
reader to \cite{AT}.\\

Our aim in this work is to give another approach to prove such kind
of existence result. The key idea in our arguments is inspired from
\cite{Bo}: By using the condition \eqref{zakaria}, we are able to
consider a suitable flow line of a considerably simplified vector
field and with some properties. As a consequence, we construct a
non-negative and bounded P-S sequence for $I_{K,\mu}$ under the
threshold
$c^\infty:=S^{\frac{n}{2}}/[n\cdot\bigl(K_\infty\bigr)^{(n-2)/2}].$
Finally, using the compactness result given in Proposition
\ref{propPS} below, we obtain a solution $u$ to problem
\eqref{problem} with\begin{equation}\label{level2}\varpi\leq
I_{K,\mu}(u)< c^\infty.\end{equation}Our arguments enable us to
regain the existence result \cite[Corollary 4.1]{Li1}. Namely, we
have the following corollary:
\begin{corollary}\label{corolla}{\it\,\,\,Let $n\geq 3.$ Assume that the hypotheses of Theorem \ref{Lions} are satisfied. Then the problem
\eqref{problem} has a solution u satisfying
\begin{equation}\label{level1}I_{K,\mu}(u)=\varpi.\end{equation}}
\end{corollary}
{\bf Remark 1.1}\,\,\,\,\eqref{level2} and \eqref{level1} enable us
to ask: What about uniqueness of critical value or, more precisely,
of the solution to \eqref{problem} in the
region $[\varpi,\,c^\infty[?$\\

As an immediate extension, we have the following existence result
which concerns the case $\mu=\mu_1(\Omega).$ Let $e_1$ denotes the
eigenfunction of $(-\Delta)$ corresponding to $\mu_1(\Omega)$ with
$e_1> 0$ and $\|e_1\|=1.$
\begin{theorem}\label{theorem}{\it\,Let $n\geq 3.$ Assume that the
hypotheses of Theorem \ref{Lions} are satisfied with $\mu
=\mu_1(\Omega).$ If $$\int_\Omega K(x)e_1^{q+1}< 0,$$ then the
problem \eqref{problem} has a solution $u$ satisfying $\varpi\leq
I_{K,\,\mu_1(\Omega)}(u)< c^\infty$ Moreover, $u$ can be chosen such
that $I_{K,\,\mu_1(\Omega)}(u)=\varpi.$}\end{theorem}

We finish this section by giving an example of function $K(x)$
dealing with the hypotheses of Theorems \ref{Lions} and
\ref{theorem}.\\\\{\bf Example 1.1}\,\,\,\, Let $y_0\in \Omega$ and
denote by $2d_0:=\mathrm{dist}(y_0,\,\partial \Omega).$ We define a
function $K:\,\bar{\Omega}\rightarrow \mathbb{R}$ by
\begin{equation*}
K(x)= -\bigl(1-\theta(|x-y_0|)\bigr)+\e_0\cdot
\theta(|x-y_0|)\bigl(d_0^\beta-\eta|x-y_0|^\beta\bigr),
\end{equation*}
where $0< \e_0< d_0,\,\,0< \eta \leq 1\,\,\text{and}\,\,\beta\geq 2$
are three fixed constants and $\theta$ is a non-increasing cut-off
function with $0\leq \theta\leq 1,$ $\theta(t)=1$ if $\,0\leq t\leq
d_0-\e_0$ and $\theta(t)=0$ if $t\geq d_0.$ A straightforward
calculation shows that the function $K(x)$
satisfies$$K_\infty=K(y_0)=\e_0 d_0^\beta> 0\quad \mathrm{and}\quad
\int_\Omega K(x)e_1^{q+1}< 0\quad\mathrm{for}\,\, \e_0\,\,\mathrm{
small\, enough}.$$On the other hand, Lions \cite[Remark 4.7]{Li1}
showed that, for $n\geq 5,$ the condition \eqref{zakaria} is
satisfied provided that $(n-2)^2\bar{c}_2\Delta
K(y_0)/\bigl(2nK(y_0)\bigr)> -\mu \bar{c}_3,$ where $\bar{c}_2$ and
$\bar{c}_3$ are two positive constants depending only on $n.$ This
means that$$\eta\frac{(n-2)^2\bar{c}_2}{d_0^2}< \mu
\bar{c}_3\quad\mathrm{for}\quad\beta=2\quad\mathrm{and}\quad \mu >
0\quad\mathrm{for}\quad\beta> 2.$$The argument used in \cite[Remark
4.7]{Li1} is still valid to show that if $n=4,$ then \eqref{zakaria}
is satisfied for any $0< \mu \leq \mu_1(\Omega).$ If $n=3,$ we
estimate the quantity $L_{K,\mu}$ by considering the test function
$u_\e$ defined by
$$u_\e(x)=\frac{\cos\bigl(\frac{\pi
|x-y_0|}{4d_0}\bigr)}{[\e+(\frac{|x-y_0|}{2d_0})^2]^{1/2}}  \quad
\mathrm {on}\,\,B(y_0,\,2d_0)\,\,\,\mathrm{and}\,\,\,
              u_\e(x)=0 \quad \mathrm{on}\,\,\Omega \setminus
B(y_0,\,2d_0),
$$where $\e> 0$ is a constant small enough, and we use
similar computations as that given in the proof of \cite[Lemma
1.3]{BN} in order to show that \eqref{zakaria} is satisfied if $\mu>
\pi^2/16d_0^2.$ This last condition is significant if, for example,
$\Omega =B(y_0,\,2d_0).$
\section{Proof of the results}\label{sec2}
\def\theequation{2.\arabic{equation}}\makeatother
\setcounter{equation}{0} {\it Proof of Theorem \ref{Lions}.}\,\,\,
To prove Theorem \ref{Lions}, we need the following result:
\begin{proposition}\label{propPS}\,\,Let $K(x)\in
C(\bar{\Omega})$ satisfying $K_\infty> 0$ and let $0< \mu\leq
\mu_1(\Omega).$ Let $c< c^\infty$ be a fixed constant. Then any
non-negative and bounded P-S sequence for $I_{K,\mu}$ at $c$ is
relatively compact.
\end{proposition}The proof is an adaptation of the arguments used to prove \cite[Lemma 1.2]{BN} and  \cite[Lemma 1]{B}. We include it in the Appendix for the
reader's convenience.\\

Let, for any $p\geq 1,$ $M_p$ denotes the following open subset of
$H_0^1(\Omega):$\begin{equation*}\label{PS9}
\begin{aligned} M_p:=\Big\{u\in
H_0^1(\Omega):\,\|u\|> (p+1)^{-1}\,\, \mathrm{and}\,\,\int_{\Omega}
Ku^{\frac{2n}{n-2}}> (p\cdot
c_\infty)^{\frac{n}{2-n}}\bigl(\int_{\Omega} |\nabla u|^2-\mu
\int_{\Omega}
u^2\bigr)^{\frac{n}{n-2}}\Big\},\end{aligned}\end{equation*} where
$c_\infty:=\bigl(n\cdot c^\infty\bigr)^{2/n}.$ The fact that
$K_\infty> 0$ and the condition \eqref{zakaria} assert that $M_p$ is
non-empty. We define the functional $J_{K,\mu}:\,M_5\rightarrow
\mathbb{R}$ by
\begin{equation}\label{prooo1}
J_{K,\mu}(u)=\frac{\int_{\Omega}|\nabla u|^2- \mu \int_{\Omega} u^2}
{\Bigl( \int_{\Omega}K|u|^{\frac{2n}{n-2}} \Bigr)^{\frac{n-2}{n}}}.
\end{equation}Let $\Sigma:=\{u\in H_0^1(\Omega):\,\|u\|=1\}$ and let $\bar{u}_0\in \bigl(M_1\cap \Sigma\bigr)$ be fixed with $\bar{u}_0\geq 0.$ Let $\tau$ be a smooth non-negative
cut-off function such that $\tau =1$ in $M_2$ and that $\tau =0$ in
$H_0^1(\Omega)\setminus M_4.$ Finally, consider the following Cauchy
problem
\begin{equation}\label{flow}
\begin{aligned}
&\frac{\partial\eta}{\partial s}(s)=W\bigl(\eta(s)\bigr),\\&
\,\,\,\,\,\eta(0)=\bar{u}_0\in M_1,
\end{aligned}\end{equation}where $W\bigl(u\bigr)$ denotes the following locally Lipschitz vector field defined by\begin{equation*}\label{vector}W\bigl(u\bigr)=\left\{
\begin{array}{ll}
\tau\bigl(u\bigr)\bigl[-\partial
J_{K,\mu}(u)\bigr] & \qquad\hbox{if} \,\,u\in M_4,\\
0 & \qquad\hbox{if}\,\,u\in H_0^1(\Omega)\setminus M_4.
\end{array}
\right.\end{equation*} Let $[0,\,T)$ denotes the positive maximal
interval defining the solution $\eta(s)$ of \eqref{flow}. We will
prove some facts satisfied by the flow line $\eta(s).$ We claim
that\begin{equation}\label{proo1}T=+\infty\quad\mathrm{and}\quad\eta(s)\in
M_1,\quad \forall\,\,s\geq 0.
\end{equation}Set $\bar{s}:=\text{sup}\{0\leq s< T:\,\,\eta(t)\in M_2,\,\,\forall\,\,0\leq t\leq
s\}.$ The continuity of the function $s\mapsto \eta(s)$ implies that
$\bar{s}> 0.$ Thus we get$$ W(\eta(s))=-\partial
J_{K,\mu}(\eta(s)),\quad \forall\,\,0\leq s< \bar{s}.$$ This,
together with \eqref{flow}, implies that $ J_{K,\mu}(\eta(s))\leq
J_{K,\mu}(\bar{u}_0)<c_\infty(K_\infty)$ for any $0\leq s< \bar{s}.$
Thus we obtain\begin{equation}\label{proo2} \eta(s)\in M_1,\quad
\forall\,\,0\leq s< \bar{s}.\end{equation} Using, again, the
continuity of the function $s\mapsto \eta(s),$ the definition of
$\bar{s},$ \eqref{proo2} and the fact that $\bar{M_1}\subset M_2$ we
derive that $\bar{s}=T.$ In particular, \eqref{flow}
becomes\begin{eqnarray}
&\qquad\quad\qquad\qquad\,\,\,\frac{\partial\eta}{\partial
s}(s)=-\partial J_{K,\mu}(\eta(s)),\quad \forall\,\,T> s\geq
0,\label{proo3}
\\& \eta(0)=\bar{u}_0\in M_1. \label{proo4}
\end{eqnarray}Using the fact that the functional $J_{K,\mu}$ is
homogenous we derive from \eqref{proo3} that
\begin{equation}\label{proo5}\bigl\langle-\partial
J_{K,\mu}(\eta(s)),\,\eta(s)\bigr\rangle=0,\quad \forall\,\,T>s\geq
0.
\end{equation}
Combining \eqref{proo3}--\eqref{proo5} we get
\begin{equation}\label{proo6}
\bigl\|\eta(s)\bigr\|=\bigl\|\bar{u}_0\bigr\|=1,\quad \forall\,\,T>
s\geq 0.
\end{equation} On
the other hand, by using the fact that $\mu< \mu_1(\Omega)$ we
derive the existence of a constant $c_0> 0$ such that
\begin{equation}\label{proo69}
\int_{\Omega}|\nabla u|^2- \mu \int_{\Omega} u^2\geq c_0,\quad
\forall\,\,u\in \Sigma.
\end{equation}The expression of $-\partial J_{K,\mu}$ and \eqref{proo69}
imply that $\|\partial J_{K,\mu}\|$ is bounded on $M_4\cap \Sigma.$
Thus $W$ is bounded on $\Sigma.$ In particular, we derive from
\eqref{proo6} that$$T=+\infty.$$This finishes the proof of the claim
\eqref{proo1}. On the other hand, \eqref{inf}, \eqref{prooo1},
\eqref{proo1}, \eqref{proo69}, together with Sobolev's inequality
and the fact that $\mathrm{sup}_{\bar{\Omega}}(K)>0,$ imply that
\begin{equation}\label{proo111}
J_{K,\mu}(\eta(s))\geq L_{K,\mu}> 0,\quad \forall\,\,s\geq 0.
\end{equation}\eqref{proo1}, \eqref{proo111} and the fact that $J_{K,\mu}(\eta(s))$ is a
non-increasing function imply that
\begin{equation}\label{proo7}L_{K,\mu}\leq \lim_{s\rightarrow
+\infty}J_{K,\mu}(\eta(s))=c<c_\infty.
\end{equation}Combining \eqref{proo3} and \eqref{proo7} we obtain
$$\int_0^{+\infty}\bigl\|\partial
J_{K,\mu}(\eta(s))\bigr\|^2\text{d}s< +\infty.$$ In particular, we
derive the existence of a sequence $(s_k),\,s_k\rightarrow +\infty,$
such that
\begin{equation}\label{proo8}\partial J_{K,\mu}(\eta(s_k))\rightarrow 0\quad
\mathrm{strongly\,\,in}\,\,H^{-1}(\Omega).
\end{equation}(In fact, by using similar arguments as that given in the proof
of \cite[Lemma A1]{BC} we can prove that $\lim_{s\rightarrow
+\infty}\|\partial J_{K,\mu}(\eta(s))\bigr\|=0$). Finally, up to
minor modifications as that given in \cite{BON} we can suppose that
\begin{equation}\label{proo9}
\eta(s_k)\geq 0,\quad \forall\,\,k\geq 0,
\end{equation}(more details will be given in a next new version).\\

Now, to prove Theorem \ref{Lions} we need to construct a sequence
$(u_k)_k$ satisfying the hypotheses of Proposition \ref{propPS}. For
this, we set, for any $k\geq 0,$
\begin{equation}\label{proo13}
u_k:=\beta_k^1\cdot\eta(s_k),
\end{equation}where $\beta_k^1:=J_{K,\mu}^{n/4}\bigl(\eta(s_k)\bigr)/\bigl(\int_{\Omega}|\nabla
\eta(s_k)|^2- \mu \int_{\Omega} |\eta(s_k)|^2\bigr)^{1/2}.$ Denoting
$$\beta_k^2:=2J_{K,\mu}^{(4-n)/4}\bigl(\eta(s_k)\bigr)/\bigl(\int_{\Omega}|\nabla
\eta(s_k)|^2- \mu \int_{\Omega} |\eta(s_k)|^2\bigr)^{1/2}.$$ A
direct calculation shows that
\begin{equation}\label{proo16}
\partial J_{K,\mu}(\eta(s_k))=\beta_k^2\cdot\partial I_{K,\mu}(u_k)\quad \mathrm{and}\quad
I_{K,\mu}(u_k)=\frac{1}{n}J_{K,\mu}^{\frac{n}{2}}(\eta(s_k)),\quad
\forall\,\,k\geq 0.
\end{equation}On the other hand,  we derive from
\eqref{proo1} and \eqref{proo6}--\eqref{proo111} the existence of
two constants $\widetilde{c}_2,\,\widetilde{c}_3> 0$ such that
\begin{equation}\label{proo14}
\widetilde{c}_2\leq \beta_k^1,\,\beta_k^2\leq \widetilde{c}_3,\quad
\forall\,\,k\geq 0.\end{equation}Combining \eqref{proo6} and
\eqref{proo7}--\eqref{proo14} we get
\begin{eqnarray}
&(u_k)_k\,\,\mathrm{is\,\,a\,\,non-negative\,\,and\,\,bounded\,\,sequence\,\,in}\,\,H_0^1(\Omega),
\nonumber\\&\varpi\leq \lim_{s\rightarrow
+\infty}I_{K,\mu}(u_k)=\frac{1}{n}c^{\frac{n}{2}}<c^\infty,
\label{proo17}\\&
\partial I_{K,\mu}(u_k)\rightarrow 0\quad
\text{strongly\,\,in}\,\,H^{-1}(\Omega). \nonumber
\end{eqnarray}These mean
that the sequence $(u_k)$ satisfies the hypotheses of Proposition
\ref{propPS}, and then, up to a subsequence, $(u_k)_k$ converges
strongly in $H_0^1(\Omega)$ to a critical point $u$ of $I_{K,\mu}$
with $u\geq 0,\,\|u\|\neq 0$ and $\varpi\leq I_{K,\mu}(u)<
c^\infty.$ It follows from the regularity theory for this kind of
equation \eqref{problem}; see, e.g., \cite[Lemma 1.5]{BN} and
\cite[Chapter 9]{B1}, that $u \in C^2(\bar{\Omega}).$ Therefore the
strong maximum principle shows that $u> 0.$ This finishes the proof
of Theorem \ref{Lions}.\\\\ {\it Proof of Corollary \ref{corolla}.}
\,\,\,Let $(w_k)$ be a non-negative minimizing sequence of the
problem \eqref{minimizing} satisfying
\begin{equation}\label{minim1}
\varpi \leq
\frac{1}{n}\bigl(\|w_k\|^2-\mu\|w_k\|_2^2\bigr)^{\frac{n}{2}}<
c^\infty,\quad \forall\,\,k.
\end{equation}
From \eqref{minim1} we can repeat the proof of Theorem \ref{Lions}
by using $w_k/\|w_k\|$ instead of $\bar{u}_0$. Therefore, from
\eqref{proo13}, \eqref{proo14} and \eqref{proo17} we obtain a
sequence $(u_k)$ of solutions for the problem \eqref{problem}
satisfying
\begin{eqnarray}
&\varpi \leq I_{K,\mu}(u_k)=\frac{1}{n}c_k^{\frac{n}{2}}\leq
\frac{1}{n}J_{K,\mu}^{\frac{n}{2}}(w_k)=
\frac{1}{n}\bigl(\|w_k\|^2-\mu\|w_k\|_2^2\bigr)^{\frac{n}{2}},
\nonumber\\& c_2\leq \|u_k\|\leq c_3. \nonumber
\end{eqnarray}
Thus $(u_k)$ is a positive and bounded P-S sequence for $I_{K,\mu}$
at $\varpi.$ This, together with Proposition \ref{propPS} and the
rest of the proof of Theorem \ref{Lions}, implies that the problem
\eqref{problem} has a solution at the level $\varpi$ for
$I_{K,\mu}.$
This finishes the proof of Corollary \ref{corolla}.\\\\
{\it Proof of Theorem \ref{theorem}.}\,\,\, To get the claims of
Theorem \ref{theorem}, it is sufficient to prove the next lemma:
\begin{lemma}\label{lemma}{\it\,
If\begin{equation}\label{theo00}\int_\Omega K(x)e_1^{q+1}<
0,\end{equation}then there exists a constant $c> 0$ such that, for
any $u\in H_0^1(\Omega)$ satisfying $\|u\|=1$ and $\int_\Omega
K(x)|u|^{q+1}\geq 0,$ we have
\begin{equation}\label{theo0}\|u\|^2-\mu_1(\Omega)\|u\|_2^2\geq c.\end{equation}}
\end{lemma}Indeed, by using Lemma \ref{lemma}
instead of the condition $\mu< \mu_1(\Omega),$ \eqref{proo69}
remains valid in $\Sigma \cap\{u\in H_0^1(\Omega):\, \int_\Omega
K(x)|u|^{q+1}\geq 0\}$ with a uniform constant $c_0.$ Now, we follow
the proof of Theorem \ref{Lions} and Corollary \ref{corolla} step by
step in order to prove that \eqref{problem} has a solution $u$ with
$\varpi\leq I_{K,\,\mu_1(\Omega)}(u)< c^\infty$ and
$I_{K,\,\mu_1(\Omega)}(u)=\varpi,$ respectively.\\\\{\it Proof of
Lemma \ref{lemma}}\,\, Arguing by contradiction, assuming that there
exists a sequence $(u_k)$ in $H_0^1(\Omega)$
satisfying\begin{eqnarray}
&\|u_k\|=1\quad \mathrm{and}\quad\int_\Omega K(x)|u_k|^{q+1}\geq
0,\quad \forall\,\,k, \label{theo1111}\\&\lim_{k\rightarrow
+\infty}\|u_k\|^2-\mu_1(\Omega)\|u_k\|_2^2=0.\label{theo1}
\end{eqnarray}Let, for
every $k,$
\begin{equation}\label{theo2}u_k=\alpha_ke_1+v_k\end{equation} be
the decomposition of $u_k$ in the Hilbert space
$\bigl(L^2(\Omega),\,\|\|_2\bigr).$ By using the Hilbert basis
$(e_s)_{s\geq 1}$ of $\bigl(L^2(\Omega),\,\|\|_2\bigr)$ defined
by$$e_s\in H_0^1(\Omega)\quad \mathrm{and}\quad
-\Delta(e_s)=\mu_s(\Omega)e_s,\quad \forall\,\,s\geq 1,$$
(\text{see,\,e.g.,\,\,\cite[Theorem 9.31]{B1}}), we derive
that\begin{equation}\label{theo3}\|v_k\|^2-\mu_1(\Omega)\|v_k\|_2^2\geq
\mathrm{inf}_{s\geq
2}\bigl(\frac{\mu_s(\Omega)-\mu_1(\Omega)}{\mu_s(\Omega)}\bigr)\|v_k\|^2,\quad
\forall\,\,k.\end{equation}Combining \eqref{theo1}--\eqref{theo3} we
obtain, up to a subsequence,\begin{equation*}\label{theo4}
\begin{aligned}0\leq
\mathrm{inf}_{s\geq
2}\bigl(\frac{\mu_s(\Omega)-\mu_1(\Omega)}{\mu_s(\Omega)}\bigr)\lim_{k\rightarrow
+\infty}\|v_k\|^2&\leq \lim_{k\rightarrow
+\infty}\|v_k\|^2-\mu_1(\Omega)\|v_k\|_2^2\\&=\lim_{k\rightarrow
+\infty}\|u_k\|^2-\mu_1(\Omega)\|u_k\|_2^2=0.\end{aligned}\end{equation*}
This, together with the fact that $\mathrm{inf}_{s\geq
2}\bigl[(\mu_s(\Omega)-\mu_1(\Omega))/\mu_s(\Omega)\bigr]\neq 0,$
implies that\begin{equation}\label{theo5}\lim_{k\rightarrow
+\infty}\|v_k\|^2=0.
\end{equation}By combining \eqref{theo1111}, \eqref{theo2} and \eqref{theo5} we
derive that$$u_k\rightarrow \pm e_1\quad
\text{strongly\,\,in}\,\,H_0^1(\Omega).$$This, together with the
continuity of the injection $H_0^1(\Omega) \subset
L^{\frac{2n}{n-2}}(\Omega)$ and \eqref{theo1111}, implies that
$$0\leq \lim_{k\rightarrow +\infty}\int_\Omega
K(x)|u_k|^{q+1}=\int_\Omega K(x)e_1^{q+1},$$ which contradicts
\eqref{theo00}. Thus the claim \eqref{theo0} follows.
\section{Appendix}\label{s3}
\def\theequation{3.\arabic{equation}}\makeatother
\setcounter{equation}{0} {\it Proof of Proposition
\ref{propPS}.}\,\,\,Let $(u_k)_k$ be a non-negative and bounded
sequence in $H_0^1(\Omega)$ satisfying
\begin{equation}\label{PS1}
I_{K,\mu}\bigl(u_k\bigr)\rightarrow c\,\,\, \mathrm{and}\,\,\,
\partial I_{K,\mu}\bigl(u_k\bigr)\rightarrow 0\,\,\, \mathrm{in}\,\,\,
H^{-1}(\Omega)\,\,\, \mathrm{with}\,\,\, c< c^\infty.
\end{equation}Since $(u_k)_k$ is bounded in
$H_0^1(\Omega),$ then there exists $u\in H_0^1(\Omega)$ such that,
up to a subsequence still denoted by $(u_k)_k,$
\begin{equation}\label{PS2}u_k\rightharpoonup u\quad \mathrm{weakly\,\,in}\quad
H_0^1(\Omega).
\end{equation}Thus, due to the fact that the injection $H_0^1(\Omega) \subset
L^2(\Omega)$ is compact, we get, up to a subsequence,
\begin{equation}\label{PS3}u_k\rightarrow u\quad \mathrm{strongly\,\,in}\quad
L^2(\Omega).
\end{equation}In particular, we derive from \eqref{PS3} that, up to
a subsequence,
\begin{equation}\label{PS4}u_k\rightarrow u\quad \mathrm{a.\,e.\,\,on}\quad
\Omega.
\end{equation}This, together with the fact that $(u_k)_k$ is bounded in
$L^{\frac{2n}{n-2}},$ implies that, passing to a further
subsequence,
\begin{equation}\label{PS5}u_k^{\frac{n+2}{n-2}}\rightharpoonup u^{\frac{n+2}{n-2}}\quad
\mathrm{weakly\,\,in}\quad L^{\frac{2n}{n+2}}\quad
(\text{see,\,\,e.g.,\,\,\cite[Exercise 4.16]{B1}}).
\end{equation}Combining \eqref{PS1}--\eqref{PS3} and \eqref{PS5} we
derive that
\begin{equation}\label{PS6}\int_{\Omega} |\nabla u|^2
- \int_{\Omega} Ku^{\frac{2n}{n-2}}- \mu \int_{\Omega} u^2=0.
\end{equation}In particular, we obtain
\begin{equation}\label{nonnegative}I_{K,\mu}(u)=\frac{1}{n}\bigl(\int_{\Omega} |\nabla u|^2
- \mu \int_{\Omega} u^2\bigr)\geq 0.
\end{equation}
Let $w_k:=u_k-u$ for any $k.$ We deduce from the Brezis--Lieb's
result \cite[Theorem 1]{BL} and \eqref{PS4}
that\begin{equation}\label{PS8}\int_{\Omega} K(x)|
u_k|^{\frac{2n}{n-2}}\text{d}x =\int_{\Omega} K(x)|
w_k|^{\frac{2n}{n-2}}\text{d}x +\int_{\Omega} K(x)|
u|^{\frac{2n}{n-2}}\text{d}x +o(1).
\end{equation}By combining \eqref{PS1}--\eqref{PS3}, \eqref{PS6} and
\eqref{PS8} we get
\begin{equation}\label{PS9}
\begin{aligned}o(1)=\bigl\langle\partial
I_{K,\mu}(u_k),\,u_k\bigr\rangle=\int_{\Omega} |\nabla u_k|^2 -
\int_{\Omega} Ku_k^{\frac{2n}{n-2}}- \mu \int_{\Omega}
u_k^2=\int_{\Omega} |\nabla w_k|^2 - \int_{\Omega}
K|w_k|^{\frac{2n}{n-2}}+o(1).
\end{aligned}\end{equation}This, together with \eqref{PS0}, implies that\begin{equation}\label{PS10}I_{K,\mu}(w_k)=\frac{1}{n}\int_{\Omega}
|\nabla w_k|^2+o(1).
\end{equation}On the other hand, by combining \eqref{PS0}, \eqref{PS2}, \eqref{PS3}, \eqref{nonnegative} and \eqref{PS8} we get
\begin{equation}\label{PS11}I_{K,\mu}(u_k)=I_{K,\mu}(w_k)+I_{K,\mu}(u)+o(1)\geq
I_{K,\mu}(w_k)+o(1).
\end{equation}We deduce from \eqref{PS1}, \eqref{PS10} and \eqref{PS11} that, for $\e>0$ a constant small enough and $k_0$ large enough,
\begin{equation}\label{PS12}\int_{\Omega}
|\nabla w_k|^2\leq
\frac{1}{\bigl(K_\infty\bigr)^{\frac{n-2}{2}}}(S-\e)^{\frac{n}{2}},\quad
\forall\,\, k\geq k_0.
\end{equation}Finally, by using Sobolev's inequality and the
fact that $K_\infty> 0$ we obtain
\begin{equation}\label{PS13}
\int_{\Omega} K|w_k|^{\frac{2n}{n-2}}\leq K_\infty\cdot
S^{-\frac{n}{n-2}}\bigl(\int_{\Omega} |\nabla
w_k|^2\bigr)^{\frac{n}{n-2}},\quad \forall\,\, k.
\end{equation}Combining \eqref{PS9}, \eqref{PS12} and \eqref{PS13},
we get
\begin{equation*}\label{PS14}\int_{\Omega}
|\nabla w_k|^2\leq
\bigl(\frac{S-\e}{S}\bigr)^{\frac{n}{n-2}}\int_{\Omega} |\nabla
w_k|^2+o(1),
\end{equation*}which implies that $\lim_{k\rightarrow
+\infty}\int_{\Omega} |\nabla w_k|^2= 0,$ and then the claim of
Proposition \ref{propPS} follows.\\\\\\
{\bf Acknowledgement.} After the accomplishment of this work in
August 2020, it was of great importance for me to find a reference
that can confirm my thought: that is why I would like to thank
Professor L. Nirenberg for his work that has gave me this
confirmation.


\begin{thebibliography}{99}
\bibitem{AT}
S. Alama and G. Tarantello, {\it On semilinear elliptic equations
with indefinite nonlinearities}, Calc. Var. Partial Differential
Equations {\bf 1} (1993), no. 4, 439--475.

\bibitem{Au}
T. Aubin, {\it Some nonlinear problems in Riemannian geometry},
Springer-Verlag, Berlin, 1998.

\bibitem{BC}
A. Bahri and J. M. Coron, {\it The scalar-curvature problem on the
standard three-dimensional sphere}, J. Funct. Anal. {\bf 95} (1991),
106--172.

\bibitem{Bo}
Z. Boucheche, {\it Existence result for an elliptic equation
involving critical exponent in three-dimensional domains}, Complex
Var. Elliptic Equ. {\bf 64} (2019), no. 4, 649--675.

\bibitem{B}
H. Brezis, {\it Elliptic equations with limiting Sobolev exponents-
the impact of topology}, Comm. Pure Appl. Math. {\bf 39} (1986), no.
S1, S17--S39.

\bibitem{B1}
H. Brezis, {\it Functional analysis, Sobolev spaces and partial
differential equations}, Springer-Verlag, New York, 2011.

\bibitem{BL}
H. Brezis and E. Lieb, {\it A relation between pointwise convergence
of functions and convergence of functionals}, Proc. Amer. Math. Soc.
{\bf 88} (1983), no. 3, 486--490.

\bibitem{BN}
H. Brezis and L. Nirenberg, {\it Positive solutions of nonlinear
elliptic equations involving critical Sobolev exponents}, Comm. Pure
Appl. Math. {\bf 36} (1983), no. 4, 437--477.

\bibitem{Li}
P. L. Lions, {\it The concentration compactness principle in the
calculus of variations. The limit case, part 1}, Rev. Mat. Iberoam.
{\bf 1} (1985), no. 1, 145--201.

\bibitem{Li1}
P. L. Lions, {\it The concentration compactness principle in the
calculus of variations. The limit case, part 2}, Rev. Mat. Iberoam.
{\bf 1} (1985), no. 2, 45--121.

\end{thebibliography}
\end{document}